\documentclass[a4paper,11pt,reqno]{amsart}
\usepackage[utf8]{inputenc}
\usepackage{amsfonts}
\usepackage{amsmath}
\usepackage{amssymb}
\usepackage{amsthm}
\usepackage{verbatim}
\usepackage{color}
\usepackage{hyperref}
\usepackage{epstopdf}

\usepackage{ mathrsfs }

\newtheorem{theorem}{Theorem}[section]

\newtheorem{lemma}[theorem]{Lemma}
\newtheorem{conjecture}[theorem]{Conjecture}
\newtheorem{proposition}[theorem]{Proposition}

\theoremstyle{remark}
\newtheorem{remark}[theorem]{Remark}

\numberwithin{equation}{section}

\renewcommand{\Re}{\textup{Re}}

\renewcommand{\H}{\mathcal{H}_{2g+1}}

\title[Mollified moments of quadratic Dirichlet $L$-functions over function fields]{Mollified moments of quadratic Dirichlet $L$-functions over function fields}

\author{Julio C. Andrade}
\address{Department of Mathematics, University of Exeter, Exeter, EX4 4QF, United Kingdom}
\email{j.c.andrade@exeter.ac.uk}

\author{Christopher G. Best}
\address{Department of Mathematics, University of Exeter, Exeter, EX4 4QF, United Kingdom}
\email{cgb212@exeter.ac.uk}

\date{\today}

\subjclass[2010]{Primary 11M38; Secondary 11M06, 11M50}
\keywords{mollified moments, non-vanishing, quadratic Dirichlet $L$-functions, function fields}

\begin{document}

\begin{abstract}
We compute asymptotic formulae for the mollified first and second moments for the family of quadratic Dirichlet $L$-functions in the function field setting. As an application, we obtain non-vanishing results for the derivatives of the completed $L$-functions $\Lambda(s,\chi_D)$ at the central point $s=1/2$. In particular, we show that the proportion of $\Lambda^{(2k)}(\frac{1}{2},\chi_D) \neq 0$ is $1+O(k^{-2})$ as $k \to \infty$.
\end{abstract}

\maketitle

\section{Introduction}

The technique of mollifying various $L$-functions has proven to be extremely useful for obtaining information about their zeros. For instance, in the case of the Riemann zeta function, Levinson \cite{Levinson1974} showed how one may obtain a lower bound on the proportion of the non-trivial zeros of $\zeta(s)$ lying on the critical line by evaluating the mollified second moment

\begin{equation} \label{zeta mollified moment}
    \frac{1}{T} \int_0^T |\zeta(\tfrac{1}{2}+it) M_N(\tfrac{1}{2}+it)|^2 dt.
\end{equation}
Here, the mollifier $M_N(s)$ is a Dirichlet polynomial of length $N$ which should approximate $\zeta(s)^{-1}$. Levinson was able to obtain an asymptotic formula as $T \to \infty$ for the mollified moment in \eqref{zeta mollified moment} for $N=T^{\theta}$ with $\theta<1/2$. As a consequence, he proved that at least $1/3$ of the non-trivial zeros of the zeta function satisfy the Riemann Hypothesis. Subsequently, Conrey \cite{Conrey1989} computed the mollified second moment for all $\theta<4/7$ and hence showed that the proportion of the non-trivial zeros on the critical line is at least $2/5$. By using different choices of mollifier, the proportion has been improved to $\geq 0.4105$ in \cite{Bui2011} and $\geq 0.417293$ in \cite{Pratt2020}.

Another problem in analytic number theory where the mollifier method has been successfully applied is the non-vanishing of $L$-functions at the central point. A famous conjecture of Chowla \cite{Chowla1965} states that $L(\frac{1}{2},\chi) \neq 0$ for all Dirichlet $L$-functions $L(s,\chi)$ with $\chi$ a primitive character. While this conjecture remains open, there are now numerous results in the literature on the non-vanishing of various families of $L$-functions at $s=1/2$. In the case of Dirichlet $L$-functions, it was shown by Iwaniec and Sarnak \cite{Iwaniec1999} that for at least $1/3$ of primitive characters $\chi$ to a sufficiently large modulus $q$, the $L$-function $L(s,\chi)$ does not vanish at $s=1/2$. This proportion was improved to $\geq 
0.341$ with the use of a two piece mollifier by Bui in \cite{Bui2012}.

In the case of the family of quadratic Dirichlet $L$-functions $L(s,\chi_d)$, Soundararajan \cite{Soundararajan2000} computed both the first and second mollified moments for a sufficiently short mollifier. This allowed Soundararajan to show that $L(\tfrac{1}{2},\chi_d) \neq 0$ for at least $7/8$ of quadratic characters $\chi_d$. By studying the one-level density of zeros of
the $L$-functions $L(s,\chi_d)$, {\"O}zl{\"u}k and Snyder \cite{Ozluk1999} showed, assuming GRH, that the proportion of non-vanishing at $s=1/2$ is at least $15/16$.

In general, the way to obtain stronger results on the zeros and non-vanishing of $L$-functions using the mollifier method is to increase the length of the mollifier. However, this is often a very difficult problem involving intricate analysis. The Ratios Conjecture of Conrey, Farmer and Zirnbauer \cite{Conrey2008} provides a relatively easy way to compute asymptotic formulae for mollified moments. In particular, one may write the mollified moments as certain contour integrals involving the mean values of ratios of $L$-functions and then apply the Ratios Conjecture to yield an asymptotic formula. Conrey and Snaith \cite{Conrey2007} gave multiple examples of computing mollified moments in this way. They show that the Ratios Conjecture implies formulae for the mollified moments with arbitrarily long mollifiers and that these formula agree with known results in the literature for sufficiently short mollifiers.

The method of mollifying may also be applied to derivatives of $L$-functions as well. For instance, it is known that RH implies that all the zeros of $\xi^{(k)}(s)$, where $\xi(s)$ is Riemann's $\xi$-function, lie on the line $\Re(s)=1/2$. By considering the second mollified moment of $\xi^{(k)}(s)$, Conrey \cite{Conrey1983} showed that the proportion of zeros of $\xi^{(k)}(s)$ on the critical line is $1+O(k^{-2})$ as $k \to \infty$. Importantly, this result does not require an arbitrarily long mollifier to be used. 

A similar approach to that employed by Soundararajan in \cite{Soundararajan2000} can be used to yield non-vanishing results for the derivatives of $L$-functions at the central point. In \cite{Kowalski2000}, Kowalski, Michel and VanderKam gave asymptotic formulae for the mollified first and second moments of the derivatives $\Lambda^{(k)}(\tfrac{1}{2},f)$ where $\Lambda(s,f)$ is the completed $L$-function attached to a primitive Hecke eigenform $f$ of weight 2. This is an orthogonal family of $L$-functions so half of the family have an even functional equation and the other half odd. Consequently, the strongest result one can hope to obtain is that $50 \%$ of $\Lambda^{(k)}(\tfrac{1}{2},f) \neq 0$ for any $k$. The key result of \cite{Kowalski2000} is that this proportion can be shown to tend to $50 \%$ as $k \to \infty$, regardless of the length of the mollifier. This is analogous to the result of Conrey \cite{Conrey1983} on the proportion of zeros of $\xi^{(k)}(s)$ on the critical line. Using similar techniques, Michel and VanderKam \cite{Michel2000} computed the mollified moments of the derivatives $\Lambda^{(k)}(\tfrac{1}{2},\chi)$ where now $\Lambda(s,\chi)$ is a completed Dirichlet $L$-function attached to the character $\chi$. They show that as $k \to \infty$, the proportion of non-vanishing of $\Lambda^{(k)}(\tfrac{1}{2},\chi)$ approaches $1/2$ if using the ``standard" mollifier and $2/3$ if they use a two-piece mollifier.

In this paper, we consider the mollified moments of the family of quadratic Dirichlet $L$-functions in the function field setting. To set some notation, let $\mathbb{F}_q$ be a finite field with $q$ odd and let $\H$ denote the hyperelliptic ensemble, i.e. the set of monic, square-free polynomials of degree $2g+1$ in $\mathbb{F}_q[t]$. For each $D \in \H$, let $\chi_D$ be the quadratic character to the modulus $D$ and let $L(s,\chi_D)$ be the associated Dirichlet $L$-function. Using a geometric argument, Li \cite{Li2018} showed that infinitely many of the $L$-functions $L(s,\chi_D)$ do vanish at $s=1/2$ and so the analogue of Chowla's conjecture in this setting does not hold. However, the density conjectures of Katz and Sarnak \cite{Katz1999} predict that $L(\tfrac{1}{2},\chi_D) \neq 0$ for $100\%$ of discriminants $D$. Bui and Florea \cite{Bui2018}, by computing the one-level density of zeros in the family, proved that the proportion of $L(\tfrac{1}{2},\chi_D)$ which do not vanish is greater than $94\%$. Ellenberg, Li and Shusterman \cite{Ellenberg2020} obtained an upper bound, depending on $q$, on the proportion of $L(s,\chi_D)$ which do vanish at $s=1/2$. Their bound tends to $0$ as $q \to \infty$ and thus improves on the bound given in \cite{Bui2018} for sufficiently large $q$. 

\section{Background and statement of results}

\subsection{Quadratic Dirichlet $L$-functions over function fields}

Here we will cover the necessary background on quadratic Dirichlet $L$-functions over function fields. Let $\mathbb{F}_q$ be a finite field with $q$ odd and let $\mathbb{F}_q[t]$ be the ring of polynomials over $\mathbb{F}_q$. We denote by $\mathcal{M}$ and $\mathcal{M}_n$ the sets of monic polynomials in $\mathbb{F}_q[t]$ and the monic polynomials of degree $n$, respectively. The polynomial $P$ is said to be prime if it is both monic and irreducible and we denote the set of prime polynomials by $\mathcal{P}$.

Throughout, we will write $d(f)$ for the degree of a polynomial $f \in \mathbb{F}_q[t]$. The norm of the polynomial $f$ is defined to be $|f|=q^{d(f)}$ if $f \neq 0$ and $|f|=0$ if $f=0$. For $\Re(s)>1$, the zeta function of $\mathbb{F}_q[t]$ is defined by the Dirichlet series and Euler product

\begin{equation}
    \zeta_q(s)=\sum_{f \in \mathcal{M}} \frac{1}{|f|^s}=\prod_{P \in \mathcal{P}} \bigg( 1-\frac{1}{|P|^s} \bigg)^{-1}.
\end{equation}
Since there are $q^n$ monic polynomials of degree $n$, we have that

\begin{equation}
    \zeta_q(s)=\sum_{n=0}^{\infty} \frac{q^n}{q^{ns}}=\frac{1}{1-q^{1-s}},
\end{equation}
which immediately provides a meromorphic continuation of $\zeta_q(s)$ to the complex plane with a simple pole at $s=1$.

We denote by $\mathcal{H}$ the set of monic, square-free polynomials in $\mathbb{F}_q[t]$ and by $\mathcal{H}_n$ the set of monic, square-free polynomials of degree $n$. The cardinality of $\mathcal{H}_n$ is

\begin{align}
    |\mathcal{H}_n|=\begin{cases}
        1 & n=0, \\
        q^n (1-q^{-1}) & n \geq 1,
    \end{cases}
\end{align}
as can be seen by considering the coefficient of $q^{-ns}$ in the series

\begin{equation}
    \sum_{n=0}^{\infty} \frac{|\mathcal{H}_n|}{q^{ns}}=\sum_{f \in \mathcal{H}} \frac{1}{|f|^s}=\frac{\zeta_q(s)}{\zeta_q(2s)}.
\end{equation}
For $D \in \H$, we define the quadratic character $\chi_D$ using the Kronecker symbol over $\mathbb{F}_q[t]$ by

\begin{equation}
    \chi_D(f)=\bigg( \frac{D}{f} \bigg).
\end{equation}
The quadratic Dirichlet $L$-function attached to the character $\chi_D$ is defined for $\Re(s)>1$ by

\begin{equation}
    L(s,\chi_D)=\sum_{f \in \mathcal{M}} \frac{\chi_D(f)}{|f|^s}=\prod_{P \in \mathcal{P}} \bigg( 1-\frac{\chi_D(P)}{|P|^s} \bigg)^{-1}.
\end{equation}
With the change of variables $u=q^{-s}$, we define $\mathcal{L}(u,\chi_D)=L(s,\chi_D)$ which is then given by the power series

\begin{equation}
    \mathcal{L}(u,\chi_D)=\sum_{f \in \mathcal{M}} \chi_D(f) u^{d(f)}=\prod_{P \in \mathcal{P}} \big( 1-\chi_D(P) u^{d(P)} \big)^{-1}.
\end{equation}
By the orthogonality relations for Dirichlet characters over $\mathbb{F}_q[t]$, we have that $\mathcal{L}(u,\chi_D)$ is in fact a polynomial in $u$ of degree $2g$ and furthermore, it satisfies the functional equation

\begin{equation}
    \mathcal{L}(u,\chi_D)=(qu^2)^g \mathcal{L} \bigg( \frac{1}{qu},\chi_D \bigg).
\end{equation}
Equivalently, we have that $L(s,\chi_D)$ is a polynomial of degree $2g$ in $q^{-s}$ and satisfies

\begin{equation}
    L(s,\chi_D)=q^{g(1-2s)} L(1-s,\chi_D).
\end{equation}
We define the completed $L$-function $\Lambda(s,\chi_D)=q^{\frac{g}{2} (2s-1)} L(s,\chi_D)$ which satisfies the symmetric functional equation

\begin{equation}
    \Lambda(s,\chi_D)=\Lambda(1-s,\chi_D).
\end{equation}

Finally, for $D \in \H$, the hyperelliptic curve $C_D$ given by the affine equation $y^2=D(x)$ is a smooth, projective and connected curve of genus $g$ over $\mathbb{F}_q$. The zeta function of the curve $C_D$

\begin{equation}
    Z_{C_D}(u)=\exp \bigg( \sum_{n=1}^{\infty} N_n(C_D) \frac{u^n}{n} \bigg),
\end{equation}
where $N_n(C_D)$ is the number of points on $C_D$ over $\mathbb{F}_{q^n}$, including the point at infinity. Weil \cite{Weil1948} showed that $Z_{C_D}(u)$ is a rational function of the form

\begin{equation}
    Z_{C_D}(u)=\frac{P_{C_D}(u)}{(1-u) (1-qu)},
\end{equation}
where $P_{C_D}(u) \in \mathbb{Z}[u]$ is a polynomial of degree $2g$. It was proven in Artin's thesis that $P_{C_D}(u)=\mathcal{L}(u,\chi_D)$. By the Riemann hypothesis for curves over finite fields, also proven by Weil \cite{Weil1948}, all of the zeros of $\mathcal{L}(u,\chi_D)$ lie on the circle $|u|=q^{-1/2}$.

\subsection{The Ratios Conjecture}

The Ratios Conjecture for the family of quadratic Dirichlet $L$-functions over function fields was formulated by Andrade and Keating in \cite{AK2014}. For simplicity, we state the conjecture below only in the case where the number of $L$-functions in the numerator and denominator is the same.

\begin{conjecture}[The Ratios Conjecture] \label{ratios conjecture}
For $|\Re(\alpha_j)|<1/4$ and $1/g \ll \Re(\beta_j)<1/4$ for $j=1,\dots,k$, we have as $g \to \infty$,

\begin{align}
    \frac{1}{|\mathcal{H}_{2g+1}|} \sum_{D\in\mathcal{H}_{2g+1}} \frac{\prod_{j=1}^k L(\tfrac{1}{2}+\alpha_j,\chi_D)}{\prod_{j=1}^k L(\tfrac{1}{2}+\gamma_j,\chi_D)}=& \sum_{\epsilon_j \in\{-1,1\}} \prod_{j=1}^k q^{g(\epsilon_j \alpha_j-\alpha_j)} YA(\epsilon_1 \alpha_1,\dots,\epsilon_k \alpha_k;\gamma)+O(q^{-\delta g})
\end{align}
for some $\delta>0$, where

\begin{equation}
    Y(\alpha;\beta)=\frac{\prod_{1\leq i\leq j\leq k} \zeta_q(1+\alpha_i+\alpha_j) \prod_{1\leq i<j\leq k} \zeta_q(1+\beta_i+\beta_j)}{\prod_{1\leq i,j\leq k} \zeta_q(1+\alpha_i+\beta_j)},
\end{equation}
and
    
\begin{align} \label{ratios conj euler product}
    A(\alpha;\beta)=& \prod_{P \in \mathcal{P}} \frac{\prod_{1\leq i\leq j\leq k} \left(1-\frac{1}{|P|^{1+\alpha_i+\alpha_j}}\right) \prod_{1\leq i<j\leq k} \left(1-\frac{1}{|P|^{1+\beta_i+\beta_j}}\right)}{\prod_{1\leq i,j\leq k} \left(1-\frac{1}{|P|^{1+\alpha_i+\beta_j}}\right)} \nonumber \\
    &\quad \times\frac{|P|}{|P|+1} \left(\frac{1}{2} \frac{\prod_{j=1}^k \left(1-\frac{1}{|P|^{1/2+\beta_j}}\right)}{\prod_{j=1}^k \left(1-\frac{1}{|P|^{1/2+\alpha_j}}\right)}+\frac{1}{2} \frac{\prod_{j=1}^k \left(1+\frac{1}{|P|^{1/2+\beta_j}}\right)}{\prod_{j=1}^k \left(1+\frac{1}{|P|^{1/2+\alpha_j}}\right)}+\frac{1}{|P|}\right).
\end{align}
\end{conjecture}

The conditions on the real parts of the shifts $\alpha_j$ and $\beta_j$ are there to ensure that the Euler products appearing in the conjectured formula are convergent. However, the conjecture may be extended to a wider ranger so long as the Euler products are convergent. As mentioned in the introduction, the above conjecture has recently been proven for certain ranges of the parameters when $k \leq 3$ by Bui, Florea and Keating in \cite{Bui2023}. Specifically, they prove the following.

\begin{theorem} \label{BFK ratios thm}
Let $0<\textup{Re}(\beta_j)<1/2$ for $1\leq j\leq k$. Denote $\alpha=\max\{|\textup{Re}(\alpha_1)|,\dots,|\textup{Re}(\alpha_k)|\}$ and $\beta=\min\{\textup{Re}(\beta_1),\dots,\textup{Re}(\beta_k)\}$. Then Conjecture \ref{ratios conjecture} holds for $1\leq k\leq 3$ with the error term $E_k$, where
    
\begin{equation}
    E_1\ll_{\varepsilon} \begin{cases} q^{-g\beta (3+2|\Re(\alpha)|)+\varepsilon g\beta} \ \text{if} \ 0\leq \Re(\alpha)<1/2 \ \text{and} \ \beta\gg g^{-1/2+\varepsilon}, \\ q^{-g\beta (3-4|\Re(\alpha)|)+\varepsilon g\beta} \ \text{if} \ -1/2<\Re(\alpha)<0 \ \text{and} \ \beta\gg g^{-1/2+\varepsilon}, \end{cases}
\end{equation}
and
    
\begin{equation}
        E_2\ll_{\varepsilon} q^{-g\beta \min\{\frac{1-4\alpha}{1+\beta},\frac{1-2\alpha}{2+\beta}\}+\varepsilon g\beta} \ \text{if} \ \alpha<1/4 \ \text{and} \ \beta\gg g^{-1/4+\varepsilon},
\end{equation}
    
\begin{equation}
    E_3\ll_{\varepsilon} q^{-g\beta \min\{\frac{1/4-4\alpha}{\beta}, \frac{1/2-4\alpha}{3+\beta}\}+\varepsilon g\beta} \ \text{if} \ \alpha<1/16 \ \text{and} \ \beta\gg g^{-1/6+\varepsilon}.
\end{equation} 
\end{theorem}

We will make use of Theorem \ref{BFK ratios thm} in the cases $k=1$ and $k=2$. We will see in the proof of our main results on the mollified moments that because Theorem \ref{BFK ratios thm} does not allow us to take $\beta \gg 1/g$ as in the Ratios Conjecture, it is the size of the bound on the error term that determines how long a mollifier we may take. For this reason, we also include the following result which improves the bound on the error $E_2$ in Theorem \ref{BFK ratios thm}.

\begin{proposition} \label{improved error prop}
With notation as in Theorem \ref{BFK ratios thm}, suppose $\alpha<1/2$ and $\beta\gg g^{-1/4+\varepsilon}$. Then, the error term $E_2$ satisfies

\begin{align}
    E_2 \ll_{\varepsilon} \begin{cases}
        q^{-2g \beta (1+2\min \{ \Re(\alpha_1), \Re(\alpha_2) \})+\varepsilon g\beta} \ &\textup{if} \ \Re(\alpha_1), \Re(\alpha_2) \geq 0, \\
        q^{-2g\beta (1-2 \max \{ |\Re(\alpha_1)|, |\Re(\alpha_2)| \})+\varepsilon g\beta} \ &\textup{if} \ \Re(\alpha_1), \Re(\alpha_2)<0, \\
        q^{-2g\beta (1+2 \min \{0, \Re(\alpha_1+\alpha_2) \})+\varepsilon g\beta} &\textup{if} \ \Re(\alpha_i) \geq 0 \ \textup{and} \ \Re(\alpha_j)<0 \ \textup{for} \ i \neq j.
    \end{cases}
\end{align}
\end{proposition}

The strategy used to prove Theorem \ref{BFK ratios thm} in \cite{Bui2023} is to first expand the $L$-functions in the denominator and write

\begin{align} \label{BFK proof start}
    & \frac{1}{|\H|} \sum_{D \in \H} \frac{\prod_{j=1}^k L(\tfrac{1}{2}+\alpha_j,\chi_D)}{\prod_{j=1}^k L(\tfrac{1}{2}+\beta_j,\chi_D)} \nonumber \\
    &\qquad =\sum_{f_1, \dots, f_k \in \mathcal{M}} \frac{\prod_{j=1}^k \mu(f_j)}{\prod_{j=1}^k |f_j|^{1/2+\beta_j}} \frac{1}{|\H|} \sum_{D \in \H} \bigg( \prod_{j=1}^k L(\tfrac{1}{2}+\alpha_j,\chi_D) \bigg) \chi_D \bigg( \prod_{j=1}^k f_j \bigg).
\end{align}
They then truncate this series and set

\begin{equation}
    S_{k, \leq X}=\sum_{f_1, \dots, f_k \in \mathcal{M}_{\leq X}} \prod_{j=1}^k \frac{\mu(f_j)}{|f_j|^{1/2+\beta_j}} \frac{1}{|\H|} \sum_{D \in \H} \bigg( \prod_{j=1}^k L(\tfrac{1}{2}+\alpha_j,\chi_D) \bigg) \chi_D \bigg( \prod_{j=1}^k f_j \bigg),
\end{equation}
where $X$ is a parameter to be chosen. Also, $S_{k, >X}$ is defined to be the sum of the terms in (\ref{BFK proof start}) where at least one of the $f_j$'s has degree larger than $X$. Then, by proving upper bounds for negative moments of these $L$-functions, it is shown in \cite{Bui2023} that

\begin{equation} \label{BFK tail bound}
    S_{k, >X} \ll_{\varepsilon} q^{-(1-\varepsilon) X \beta},
\end{equation}
for $\beta \gg g^{-1/2k+\varepsilon}$. Also proven in \cite{Bui2023} are asymptotic formulae for the twisted, shifted moments appearing in (\ref{BFK proof start}) for $k \leq 3$. Inserting the main terms of these formulae into $S_{k, \leq X}$, extending the sums over $f_j$ to be over all $f_j \in \mathcal{M}$ (which introduces a negligible error) and then performing an Euler product computation yields the main terms as predicted by the Ratios Conjecture.

For the contribution of the error terms in the twisted moment formulae to $S_{k, \leq X}$ in the case $k=2$ or $k=3$, Bui, Florea and Keating use their overall bounds for the error and bound the sums over $f_j$ trivially. An optimum value for the parameter $X$ is then chosen to yield the bound on the error terms $E_2$ and $E_3$ in Theorem \ref{BFK ratios thm}. However, for $k=1$, they keep the error terms in the first twisted moment explicit and make use of the cancellation coming from the M{\"o}bius function in the sum over $f_1$. This leads to a better error term and a wider range for the shift parameter $\alpha$. The proof of Proposition \ref{improved error prop} is similar to the case of $k=1$ covered in \cite[Section 4.1]{Bui2023} and so we omit the details.

\subsection{The mollification}

For our choice of mollifier, we take the standard notion of a mollifier in the number field setting and translate it to the function field setting in the natural way. The mollifier is a Dirichlet polynomial whose purpose is to approximate $L(\tfrac{1}{2},\chi_D)^{-1}$ so we begin with the Dirichlet series

\begin{equation}
    \frac{1}{L(s,\chi_D)}=\sum_{f \in \mathcal{M}} \frac{\mu(f) \chi_D(f)}{|f|^s},
\end{equation}
which is convergent for $\Re(s)>1/2$. Truncating this series and multiplying by a smoothing function leads us to the mollifier

\begin{equation}
    M_y(\chi_D,P)=\sum_{\substack{f \in \mathcal{M} \\ |f| \leq y}} \frac{\mu(f) \chi_D(f)}{|f|^{1/2}} P\left(\frac{\log(y/|f|)}{\log y}\right),
\end{equation}
where $P(x)$ is a polynomial satisfying $P(0)=0$ and $y=(q^{2g})^{\theta}$ for some $\theta>0$. Provided that the mollifier is not too long, i.e. $\theta$ is sufficiently small, we will prove asymptotic formula for the first and second mollified moments

\begin{equation} \label{mollified moment 1}
    \frac{1}{|\H|} \sum_{D \in \H} L(\tfrac{1}{2},\chi_D) M_y(\chi_D,P),
\end{equation}
and

\begin{equation} \label{mollified moment 2}
    \frac{1}{|\H|} \sum_{D \in \H} |L(\tfrac{1}{2},\chi_D) M_y(\chi_D,P)|^2.
\end{equation}
Our main results are now stated below.




\begin{theorem} \label{1st mollified moment thm}
Let $Q$ be an even polynomial and $P$ a polynomial satisfying $P(0)=0$. Then for $\theta<3/2$, we have as $g \to \infty$,

\begin{align}
    & Q \left( \frac{1}{g\log q}\frac{d}{d\alpha} \right) \frac{1}{|\H|} \sum_{D\in\mathcal{H}_{2g+1}} \Lambda(\tfrac{1}{2}+\alpha,\chi_D) M_y(\chi_D,P) \bigg|_{\alpha=0}=P(1) Q(1)+\frac{1}{2\theta} P'(1) \int_0^1 Q(t) \, dt+O(1/g).
\end{align}
Assuming Conjecture \ref{ratios conjecture}, the result holds for any $\theta>0$.
\end{theorem}

\begin{theorem} \label{2nd mollified moment thm}
Let $Q_1, Q_2$ be even polynomials and $P_1, P_2$ polynomials satisfying $P_j(0)=P_j'(0)=0$ for $j=1,2$. Then for $\theta<1/2$, we have as $g \to \infty$,

\begin{align}
    & Q_1 \left( \frac{1}{g\log q} \frac{d}{d\alpha} \right) Q_2 \left( \frac{1}{g\log q} \frac{d}{d\beta} \right) \frac{1}{|\H|} \sum_{D\in\mathcal{H}_{2g+1}} \Lambda(\tfrac{1}{2}+\alpha,\chi_D) \Lambda(\tfrac{1}{2}+\beta,\chi_D) M_y(\chi_D,P_1) M_y(\chi_D,P_2)\bigg|_{\alpha=\beta=0} \nonumber \\
    & =\frac{1}{8\theta} \int_0^1\int_0^1 \bigg(\frac{1}{\theta} P_1''(r) \Tilde{Q}_1(u)-4\theta P_1(r) Q_1'(u)\bigg) \left(\frac{1}{\theta} P_2''(r) \Tilde{Q}_2(u)-4\theta P_2(r) Q_2'(u)\right) \, du \, dr \nonumber \\
    &\quad +\frac{1}{4}\left(\frac{1}{\theta} P_1'(1) \Tilde{Q}_1(1)+2P_1(1) Q_1(1)\right) \left(\frac{1}{\theta} P_2'(1) \Tilde{Q}_2(1)+2P_2(1) Q_2(1)\right)+O(1/g),
\end{align}
where
    
\begin{equation}
    \Tilde{Q}(u)=\int_0^u Q(t) \, dt.
\end{equation}
Assuming Conjecture \ref{ratios conjecture}, the result holds for any $\theta>0$.
\end{theorem}

\begin{remark}
\begin{enumerate}
    \item We only consider even polynomials $Q$ in Theorems \ref{1st mollified moment thm} and \ref{2nd mollified moment thm} since by the functional equation

    \begin{equation}
        \Lambda(s,\chi_D)=\Lambda(1-s,\chi_D),
    \end{equation}
    we have that $\Lambda^{(k)}(\tfrac{1}{2},\chi_D)=0$ if $k$ is odd.

    \item Theorem \ref{2nd mollified moment thm} is the function field analogue of Theorem 5.2 in \cite{Conrey2007} which gives a formula for the mollified second moment of quadratic Dirichlet $L$-functions assuming the Ratios Conjecture.

    \item Asymptotics for the mollified moments in \eqref{mollified moment 1} and \eqref{mollified moment 2} can be recovered by taking $Q(x)=Q_1(x)=Q_2(x)=1$ in Theorems \ref{1st mollified moment thm} and \ref{2nd mollified moment thm}. In particular, we find that

    \begin{equation}
        \frac{1}{|\H|} \sum_{D \in \H} L(\tfrac{1}{2},\chi_D) M_y(\chi_D,P) \sim P(1)+\frac{1}{2\theta} P'(1),
    \end{equation}
    and

    \begin{equation}
        \frac{1}{|\H|} \sum_{D \in \H} |L(\tfrac{1}{2},\chi_D) M_y(\chi_D,P)|^2 \sim \left( P(1)+\frac{1}{2\theta} P'(1) \right)^2+\frac{1}{24\theta^3} \int_0^1 P''(r)^2 \, dr.
    \end{equation}
\end{enumerate}
\end{remark}

\subsection{Applications to non-vanishing}

As an application of our results on the mollified moments, we are able to obtain non-vanishing results on the derivatives $\Lambda^{(2k)}(\tfrac{1}{2},\chi_D)$. For $Q(x)=\sum_{n \geq 0} a_n x^n$ an even polynomial, we define the differential operator

\begin{equation}
    \widehat{Q}=Q \left( \frac{1}{g \log q} \frac{d}{ds} \right)=\sum_{n \geq 0} \frac{a_n}{(g \log q)^n} \frac{d^n}{ds^n}.
\end{equation}
Then we define the first and second mollified moments

\begin{equation}
    \mathcal{S}_1(P,Q)=\frac{1}{|\H|} \sum_{D \in \H} \widehat{Q} \big( \Lambda(s,\chi_D) \big) (\tfrac{1}{2}) M_y(\chi_D,P),
\end{equation}
and

\begin{equation}
    \mathcal{S}_2(P,Q)=\frac{1}{|\H|} \sum_{D \in \H} |\widehat{Q} \big( \Lambda(s,\chi_D) \big) (\tfrac{1}{2}) M_y(\chi_D,P)|^2.
\end{equation}
By Theorems \ref{1st mollified moment thm} and \ref{2nd mollified moment thm}, we have

\begin{equation}
    \mathcal{S}_1(P,Q)=P(1) Q(1)+\frac{1}{2\theta} P'(1) \Tilde{Q}(1)+O(1/g),
\end{equation}
and

\begin{align}
    \mathcal{S}_2(P,Q) & =\bigg( P(1) Q(1)+\frac{1}{2\theta} P'(1) \Tilde{Q}(1) \bigg)^2 \nonumber \\
    &\qquad +\frac{1}{8\theta} \int_0^1 \int_0^1 \bigg( \frac{1}{\theta} P''(r) \Tilde{Q}(u)-4\theta P(r) Q'(u) \bigg)^2 \, du \, dr+O(1/g).
\end{align}
By applying the Cauchy-Schwarz inequality in the usual way, we have

\begin{equation}
    \frac{1}{|\H|} |\{ D \in \H: \widehat{Q} \big( \Lambda(s,\chi_D) \big) (\tfrac{1}{2}) \neq 0 \}| \geq \frac{\mathcal{S}_1(P,Q)^2}{\mathcal{S}_2(P,Q)}.
\end{equation}
Using the above two expression for $\mathcal{S}_1(P,Q)$ and $\mathcal{S}_2(P,Q)$ on the right-hand side of this inequality then gives us the following general non-vanishing result.

\begin{theorem} \label{general non-vanishing thm}
For $Q$ an even polynomial, we have that as $g \to \infty$,

\begin{equation} \label{general non-vanishing bound}
    \frac{1}{|\H|} |\{ D \in \H: \widehat{Q} \big( \Lambda(s,\chi_D) \big) (\tfrac{1}{2}) \neq 0 \}| \geq \sup_{P,\theta} \frac{1}{1+\mathcal{R}(P,Q)}+o(1),
\end{equation}
where the supremum is over all polynomials $P$ with $P(0)=P'(0)=0$ and real numbers  $0<\theta<1/2$. Also, the ratio $\mathcal{R}(P,Q)$ is given by

\begin{equation} \label{ratio def}
    \mathcal{R}(P,Q)=\frac{(2\theta)^{-1} \int_0^1 \int_0^1 \big( P''(r) \Tilde{Q}(u)-4\theta^2 P(r) Q'(u) \big)^2 \, du \, dr}{\big( 2\theta P(1) Q(1)+P'(1) \Tilde{Q}(1) \big)^2}.
\end{equation}
\end{theorem}

\begin{remark} \label{non-vanishing remark}
\begin{enumerate}
    \item The ratio $\mathcal{R}(P,Q)$ is non-negative so the non-vanishing proportion is at most 1 as expected.

    \item The expression for the ratio $\mathcal{R}(P,Q)$ above is very similar to that for the analogous ratio appearing in \cite{Kowalski2000}, see \cite[Eq. 31]{Kowalski2000}. 

    \item By the continuity of $\mathcal{R}(P,Q)$ in $P$ and $\theta$, when determining the value of the supremum above, one may take $\theta=1/2$ and allow $P$ to range over all power series $P(x)=a_2x^2+a_3x^3+\dots$ such that this series and the series for everything up to and including $P''(x)^2$ are absolutely convergent on $[0,1]$.
\end{enumerate}
\end{remark}

Finally, we examine in more detail the non-vanishing of $\Lambda^{(2k)}(\tfrac{1}{2},\chi_D)$. We denote the proportion of non-vanishing by

\begin{equation}
    p_{2k}=\liminf_{g \to \infty} \frac{|\{ D \in \H: \Lambda^{(2k)}(\tfrac{1}{2},\chi_D) \neq 0 \}|}{|\H|}.
\end{equation}
By taking $Q(x)=x^{2k}$ in Theorem \ref{general non-vanishing thm} and determining the optimal choice of $P$ to minimise $\mathcal{R}(P,x^{2k})$, we will prove the following.

\begin{theorem} \label{non-vanishing thm}
For all $k \geq 0$, we have that $p_{2k}>0$. Moreover, we have $p_{2k} \geq \pi_{2k}$, where $\pi_{2k}$ satisfies

\begin{equation}
    \pi_{2k}=1-\frac{1}{16 (2k)^2}+O(k^{-3}),
\end{equation}
as $k \to \infty$. In particular, we have

\begin{equation}
     p_0 \geq 0.875, \ p_2 \geq 0.9895, \ p_4 \geq 0.9971, \ p_6 \geq 0.9986 \ \textup{and} \ p_8 \geq 0.9991.
\end{equation}
\end{theorem}

Our lower bound of $0.875$ on the proportion $p_0$ of non-vanishing of $L(\tfrac{1}{2},\chi_D)$ is not an improvement on the results of \cite{Bui2018, Ellenberg2020} but shows a clear analogy between our results and those of Soundararajan in \cite{Soundararajan2000}. Our Theorem \ref{non-vanishing thm} is also a function field analogue of \cite[Theorem 1.2]{Kowalski2000} and \cite[Theorem 1]{Michel2000} and shows that $L$-functions $L(s,\chi_D)$ with a high order of vanishing at $s=1/2$ are rare.

\section{The mollified first moment}

In this section we will prove Theorem \ref{1st mollified moment thm}. We begin by considering the shifted mollified moment

\begin{equation} \label{1st shifted mollified def}
    \mathcal{M}(\alpha,P):=\frac{1}{|\H|} \sum_{D\in\mathcal{H}_{2g+1}} L(\tfrac{1}{2}+\alpha,\chi_D) M_y(\chi_D,P).
\end{equation}

\begin{lemma} \label{1st shifted mollified lemma}
For $\theta<3/2$, we have
    
\begin{equation}
    \mathcal{M}(\alpha,P)=\frac{1+q^{-2g\alpha}}{2} P(1)+\frac{1-q^{-2g\alpha}}{2\alpha \log y} P'(1)+O(1/g),
\end{equation}
uniformly for $\alpha\ll 1/g$. Assuming Conjecture \ref{ratios conjecture}, the result holds for any $\theta>0$.
\end{lemma}

\begin{proof}
First, we write the mollifier in the form of a contour integral by using the fact that for $|f| \leq y$ and $n \in \mathbb{N}$,

\begin{equation}
    \left( \frac{\log(y/|f|)}{\log y} \right)^n=\frac{n!}{\log^n y} \frac{1}{2\pi i} \int_{(c)} \left(\frac{y}{|f|} \right)^z \frac{dz}{z^{n+1}},
\end{equation}
where $c>0$. This can be seen by moving the contour to $-\infty$ and computing the residue of the pole at $z=0$. Therefore, by writing $P(x)=\sum_{n \geq 1} p_n x^n$, we have

\begin{equation} 
    M_y(\chi_D,P)=\sum_{n\geq 1} \frac{p_n n!}{\log^n y} \sum_{\substack{f \in \mathcal{M} \\ |f| \leq y}} \frac{\mu(f) \chi_D(f)}{|f|^{1/2}} \frac{1}{2\pi i} \int_{(c)} \left( \frac{y}{|f|} \right)^z \frac{dz}{z^{n+1}}.
\end{equation}
Now, if $|f|>y$, by moving the contour far to the right we see that the integral above vanishes. Thus, we may drop the condition that $|f| \leq y$ from the sum over $f$. Then, since $\textrm{Re}(z)>0$, the resulting series converges and we have

\begin{equation} \label{mollifier integral}
    M_y(\chi_D,P)=\sum_{n\geq 1} \frac{p_n n!}{\log^n y} \frac{1}{2\pi i} \int_{(c)} \frac{y^z}{z^{n+1}} \frac{dz}{L(\tfrac{1}{2}+z,\chi_D)}.
\end{equation}
Using this expression for the mollifier in (\ref{1st shifted mollified def}) yields

\begin{equation} \label{I(P) integral}
    \mathcal{M}(\alpha,P)=\sum_{n \geq 1} \frac{p_n n!}{\log^n y} \frac{1}{2\pi i} \int_{(c)} \frac{y^z}{z^{n+1}} \frac{1}{|\H|} \sum_{D \in \H} \frac{L(\tfrac{1}{2}+\alpha,\chi_D)}{L(\tfrac{1}{2}+z,\chi_D)} \, dz.
\end{equation}

By Theorem \ref{BFK ratios thm}, for $|\textup{Re}(\alpha)|<1/2$ and $g^{-1/2+\varepsilon} \ll c<1/2$, the mean value of the ratio of $L$-functions appearing in the integrand is given by

\begin{align} \label{1 over 1 ratio}
    \frac{1}{|\H|} \sum_{D \in \H} \frac{L(\tfrac{1}{2}+\alpha,\chi_D)}{L(\tfrac{1}{2}+z,\chi_D)} & =\frac{\zeta_q(1+2\alpha)}{\zeta_q(1+\alpha+z)} A(\alpha;z)+q^{-2g\alpha} \frac{\zeta_q(1-2\alpha)}{\zeta_q(1-\alpha+z)} A(-\alpha;z) \nonumber \\
    &\qquad +O(q^{-gc(3-4|\textup{Re}(\alpha)|)+\varepsilon cg}),
\end{align}
where $A(\alpha;z)$ is the Euler product defined in (\ref{ratios conj euler product}). As $y=q^{2g \theta}$ and $\alpha \ll 1/g$, by bounding the integral by absolute values, the contribution of the error term above to \eqref{I(P) integral} will be

\begin{equation}
    \ll_{\varepsilon} q^{cg (2\theta -3+\varepsilon)}.
\end{equation}
This error is $\ll_{\varepsilon} q^{-\varepsilon g}$ for some $\varepsilon>0$ if and only if $\theta<3/2$. On the other hand, if we assume the Ratios Conjecture, then \eqref{1 over 1 ratio} holds for $1/g \ll c<1/2$ with an error that is $O(q^{-\varepsilon g})$ uniformly. In this case, we may take $c \asymp 1/g$ and then the contribution of the error term to \eqref{I(P) integral} will be $\ll_{\varepsilon} q^{-\varepsilon g}$ for any $\theta>0$. Thus, unconditionally for $\theta<3/2$ and conditionally for any $\theta>0$, we may write

\begin{equation} \label{1st M(alpha;P) split}
    \mathcal{M}(\alpha,P)=I(\alpha,P)+q^{-2g\alpha} I(-\alpha,P)+O(q^{-\varepsilon g}),
\end{equation}
where

\begin{equation} \label{1st I def}
    I(\alpha,P)=\zeta_q(1+2\alpha) \sum_{n \geq 1} \frac{p_n n!}{\log^n y} J_{\alpha}(y),
\end{equation}
with

\begin{equation}
    J_{\alpha}(y)=\frac{1}{2\pi i} \int_{(c)} \frac{y^z}{z^{n+1}} \frac{A(\alpha;z)}{\zeta_q(1+\alpha+z)} \, dz.
\end{equation}

By moving the contour to $\textup{Re}(z)=-\delta$, where $\delta>0$ is sufficiently small so that the Euler product is absolutely convergent, we have that $J_{\alpha}(y)$ is given by the residue at $z=0$ plus the new integral along the line $\textup{Re}(z)=-\delta$. We write the residue at $z=0$ as an integral where the contour is a circle around zero and we bound the new integral trivially by $O(y^{-\delta})$. We therefore have that

\begin{equation}
    J_{\alpha}(y)=\frac{1}{2\pi i} \oint \frac{y^z}{z^{n+1}} \frac{A(\alpha;z)}{\zeta_q(1+\alpha+z)} \, dz+O(q^{-\varepsilon g}),
\end{equation}
where the contour is now a circle of radius $\asymp 1/g$ around zero. On this circular contour, for $\alpha, z \ll 1/g$, we use the Taylor expansion

\begin{equation}
    \frac{A(\alpha;z)}{\zeta_q(1+\alpha+z)}=(\alpha+z) A(0;0) \log q+O(g^{-2}).
\end{equation}
From the definition in (\ref{ratios conj euler product}), we have that $A(0;0)=1$ which gives us

\begin{equation} \label{1st J expression}
    J_{\alpha}(y)=\frac{\log q}{2\pi i} \oint \frac{y^z}{z^{n+1}} \left( (\alpha+z)+O(g^{-2}) \right) \, dz=\frac{\log q}{2\pi i} \oint \frac{y^z}{z^{n+1}} (\alpha+z) \, dz+O(g^{n-2}).
\end{equation}
Next, by computing the residue at $z=0$, we have

\begin{equation} \label{residue 1}
    \sum_{n\geq 1} \frac{p_n n!}{\log^n y} \frac{1}{2\pi i}\oint \frac{y^z}{z^{n+1}} dz=\sum_{n \geq 1} p_n=P(1)
\end{equation}
and

\begin{equation} \label{residue 2}
    \sum_{n\geq 1} \frac{p_n n!}{\log^n y} \frac{1}{2\pi i}\oint \frac{y^z}{z^n} dz=\frac{1}{\log y} \sum_{n \geq 1} np_n=\frac{1}{\log y} P'(1).
\end{equation}
Therefore, by combining \eqref{1st I def}, \eqref{1st J expression}, \eqref{residue 1} and \eqref{residue 2}, we have

\begin{equation}
    I(\alpha;P)=\zeta_q(1+2\alpha) \log q \left( \alpha P(1)+\frac{1}{\log y} P'(1)+O(g^{-2}) \right).
\end{equation}

Lastly, for $\alpha \asymp 1/g$, we have the Laurent expansion

\begin{equation}
    \zeta_q(1+2\alpha)=\frac{1}{2\alpha \log q}+O(1)=O(g),
\end{equation}
and so we have

\begin{equation}
    I(\alpha;P)=\frac{1}{2} P(1)+\frac{1}{2\alpha \log y} P'(1)+O(g^{-1}),
\end{equation}
uniformly on any fixed annulus $|\alpha| \asymp 1/g$. Plugging this back into (\ref{1st M(alpha;P) split}) yields

\begin{equation}
    \mathcal{M}(\alpha;P)=\frac{1+q^{-2g\alpha}}{2} P(1)+\frac{1-q^{-2g\alpha}}{2\alpha \log y} P'(1)+O(g^{-1}),
\end{equation}
with the assumption that $\alpha \asymp 1/g$. But since $\mathcal{M}(\alpha;P)$ and the main term above are holomorphic for $\alpha \ll 1/g$, the error term is also holomorphic in this region. By the maximum modulus principle, the bound on the error term also holds uniformly for $\alpha \ll 1/g$ which completes the proof.
\end{proof}

\subsection{Proof of Theorem \ref{1st mollified moment thm}}

We define

\begin{equation}
    \mathcal{N}(\alpha;P)=\frac{1}{|\H|} \sum_{D \in \H} \Lambda(\tfrac{1}{2}+\alpha,\chi_D) M_y(\chi_D,P).
\end{equation}
By definition, $\Lambda(\tfrac{1}{2}+\alpha,\chi_D)=q^{g\alpha} L(\tfrac{1}{2}+\alpha,\chi_D)$ and so by Lemma \ref{1st shifted mollified lemma}, we have

\begin{equation}
    \mathcal{N}(\alpha;P)=\frac{q^{g\alpha}+q^{-g\alpha}}{2} P(1)+\frac{q^{g\alpha}-q^{-g\alpha}}{2\alpha \log y} P'(1)+O(g^{-1}),
\end{equation}
uniformly for $\alpha \ll 1/g$. We now write $\alpha=\frac{a}{g \log q}$. Then, as $y=(q^{2g})^{\theta}$, the above can be rewritten as

\begin{equation}
    \mathcal{N} \left( \frac{a}{g \log q};P \right)=P(1) \cosh{a}+P'(1) \frac{\sinh{a}}{2a \theta}+O(g^{-1}).
\end{equation}
For an even polynomial $Q$, we have that

\begin{equation}
    Q \left( \frac{d}{da} \right) \left. \cosh{a} \right|_{a=0}=Q(1),
\end{equation}
and

\begin{equation}
    Q \left( \frac{d}{da} \right) \left. \frac{\sinh{a}}{a} \right|_{a=0}=Q \left( \frac{d}{da} \right) \left. \int_0^1 \cosh{at} \, dt \right|_{a=0}=\int_0^1 Q(t) \, dt.
\end{equation}
Thus, we have

\begin{equation}
    Q \left( \frac{d}{da} \right) \left. \mathcal{N} \left( \frac{a}{g \log q};P \right) \right|_{a=0}=P(1) Q(1)+\frac{1}{2\theta} P'(1) \int_0^1 Q(t) \, dt+O(g^{-1}).
\end{equation}
By noting that $\frac{d}{da}=\frac{1}{g \log q} \frac{d}{d\alpha}$, we see that this is precisely the statement of Theorem \ref{1st mollified moment thm}.

\section{The mollified second moment}

In this section we will prove Theorem \ref{2nd mollified moment thm} on the mollified second moment. Similarly to the first moment, we will first obtain a formula for the shifted moment

\begin{equation}
    \mathcal{M}(\alpha,\beta;P_1,P_2):=\frac{1}{|\H|} \sum_{D\in\mathcal{H}_{2g+1}} L(\tfrac{1}{2}+\alpha,\chi_D) L(\tfrac{1}{2}+\beta,\chi_D) M_y(\chi_D,P_1) M_y(\chi_D,P_2).
\end{equation}

\begin{lemma} \label{2nd shifted mollified lemma}
    For polynomials $P_1, P_2$ satisfying $P_j(0)=P_j'(0)=0$ for $j=1,2$, and for $\theta<1/2$, we have

    \begin{align}
    \mathcal{M}(\alpha,\beta;P_1,P_2) & =\frac{1}{4} \left(\frac{\alpha\beta (1-q^{-2g(\alpha+\beta)})}{\alpha+\beta}+\frac{\alpha\beta (q^{-2g\alpha}-q^{-2g\beta})}{\alpha-\beta}\right) \log y \int_0^1 P_1(r) P_2(r) dr \nonumber \\
    &\quad +\frac{1}{4} (1+q^{-2g\alpha}) (1+q^{-2g\beta}) \int_0^1 \big(P_1(r) P_2'(r) +P_1'(r) P_2(r)\big) dr \nonumber \\
    &\quad +\frac{1}{4} \left(\frac{(1+q^{-2g\alpha}) (1-q^{-2g\beta})}{\beta}+\frac{(1-q^{-2g\alpha}) (1+q^{-2g\beta})}{\alpha}\right) \frac{1}{\log y} \int_0^1 P_1'(r) P_2'(r) dr \nonumber \\
    &\quad +\frac{1}{4} \left(\frac{1-q^{-2g(\alpha+\beta)}}{\alpha+\beta}+\frac{q^{-2g\beta}-q^{-2g\alpha}}{\alpha-\beta}\right) \frac{1}{\log y} \int_0^1 \big(P_1(r) P_2''(r)+P_1''(r) P_2(r)\big) dr \nonumber \\
    &\quad +\frac{(1-q^{-2g\alpha}) (1-q^{-2g\beta})}{4\alpha\beta} \frac{1}{\log^2 y} \int_0^1 \big(P_1'(r) P_2''(r)+P_1''(r) P_2'(r)\big) dr \nonumber \\
    &\quad +\frac{1}{4\alpha\beta} \left(\frac{1-q^{-2g(\alpha+\beta)}}{\alpha+\beta}+\frac{q^{-2g\alpha}-q^{-2g\beta}}{\alpha-\beta}\right) \frac{1}{\log^3 y} \int_0^1 P_1''(r) P_2''(r) dr+O(g^{-1})
    \end{align}
    uniformly for $\alpha,\beta \ll 1/g$. Assuming Conjecture \ref{ratios conjecture}, the result holds for all $\theta>0$.
\end{lemma}

\begin{proof}
We write $P_1(x)=\sum_{m\geq 2} p_{1,m} x^m$ and $P_2(x)=\sum_{n\geq 2} p_{2,n} x^n$ and use \eqref{mollifier integral} which gives us

\begin{align} \label{I(P,P) integral}
    \mathcal{M}(\alpha,\beta;P_1,P_2) & =\sum_{m,n\geq 2} \frac{p_{1,m} m! p_{2,n} n!}{\log^{m+n} y} \frac{1}{(2\pi i)^2}\int_{(c)}\int_{(c)} \frac{y^{w+z}}{w^{m+1} z^{n+1}} \nonumber \\ 
    &\qquad \times\frac{1}{|\mathcal{H}_{2g+1}|} \sum_{D\in\mathcal{H}_{2g+1}} \frac{L(\tfrac{1}{2}+\alpha,\chi_D) L(\tfrac{1}{2}+\beta,\chi_D)}{L(\tfrac{1}{2}+w,\chi_D) L(\tfrac{1}{2}+z,\chi_D)} \, dw \, dz,
\end{align}
where $c>0$. For $\alpha,\beta \ll 1/g$ and $g^{-1/4+\varepsilon} \ll c<1/2$, Proposition \ref{improved error prop} gives us that the mean value of the ratio of $L$-functions in the integrand is given by Conjecture \ref{ratios conjecture} with an error that is $\ll_{\varepsilon} q^{-2cg+\varepsilon cg}$. By bounding the integral with absolute values, the contribution of this error to \eqref{I(P,P) integral} will be

\begin{equation}
    \ll_{\varepsilon} y^{2c} q^{-2cg+\varepsilon cg} \ll_{\varepsilon} q^{2cg(2\theta-1+\varepsilon/2)}.
\end{equation}
For $\theta<1/2$, we may bound this error by $O(q^{-\varepsilon g})$ for sufficiently small $\varepsilon>0$. Assuming Conjecture \ref{ratios conjecture}, the error in the ratio of $L$-functions is $\ll_{\varepsilon} q^{-\varepsilon g}$ and we may take $c \asymp 1/g$. This allows us to bound the error in \eqref{I(P,P) integral} by $O(q^{-\varepsilon g})$ for any $\theta>0$. In either case, by inserting the main terms of the Ratios Conjecture into \eqref{I(P,P) integral}, we write

\begin{align} \label{M(a,b) decomp}
    \mathcal{M}(\alpha,\beta;P_1,P_2) & =I(\alpha,\beta)+q^{-2g\alpha} I(-\alpha,\beta)+q^{-2g\beta} I(\alpha,-\beta)+q^{-2g(\alpha+\beta)} I(-\alpha,-\beta)+O(q^{-\varepsilon g}),
\end{align}
where

\begin{align} \label{I(a,b) def}
    I(\alpha,\beta)=\zeta_q(1+2\alpha) \zeta_q(1+2\beta) \zeta_q(1+\alpha+\beta) \sum_{m,n\geq 2} \frac{p_{1,m} m! p_{2,n} n!}{(\log y)^{m+n}} J_{\alpha,\beta}(y),
\end{align}
with

\begin{align}
    J_{\alpha,\beta}(y)& =\frac{1}{(2\pi i)^2}\int_{(c)}\int_{(c)} \frac{y^{w+z}}{w^{m+1} z^{n+1}} \frac{\zeta_q(1+w+z) A(\alpha,\beta;w,z)}{\zeta_q(1+\alpha+w) \zeta_q(1+\alpha+z) \zeta_q(1+\beta+w) \zeta_q(1+\beta+z)} \, dw \, dz.
\end{align}

To evaluate the contour integral $J_{\alpha,\beta}(y)$, we would like to move the contours to the left of 0 but the factor of $\zeta_q(1+w+z)$ with its pole at $w+z=0$ makes this somewhat tricky. To deal with this, we use the fact that for $\Re(w+z)>0$, we may write

\begin{equation}
    \frac{y^{w+z}}{(w+z)}=\int_0^y u^{w+z} \frac{du}{u}.
\end{equation}
Thus, we have

\begin{align}
    J_{\alpha,\beta}(y) & =\frac{1}{(2\pi i)^2} \int_1^y \int_{(c)} \int_{(c)} \frac{u^{w+z}} {w^{m+1} z^{n+1}} \nonumber \\
    &\qquad \times \frac{(w+z) \zeta_q(1+w+z) A(\alpha,\beta;w,z)}{\zeta_q(1+\alpha+w) \zeta_q(1+\alpha+z) \zeta_q(1+\beta+w) \zeta_q(1+\beta+z)} \, dw \, dz \frac{du}{u},
\end{align}
with the integration in $u$ only over $1 \leq u \leq y$ since if $u<1$, we can move the contours far to the right to see that the integrals in $w$ and $z$ vanish. As $(w+z) \zeta_q(1+w+z)$ is analytic at $w+z=0$, the poles of the integrand are now at $w=0$ or $z=0$ only. So, by moving the contours to $\Re(w)=\Re(z)=-\delta$ where $\delta>0$ is sufficiently small to ensure that the Euler product is absolutely convergent, we get that $J_{\alpha,\beta}(y)$ is given by the residue at $w=z=0$ plus the new integral over the lines $\Re(w)=\Re(z)=-\delta$ whose contribution may be bounded by

\begin{equation}
    \int_1^y u^{-2\delta} \frac{du}{u} \ll 1.
\end{equation}

We express the residue at $w=z=0$ as an integral with the contours being circles of radius $\asymp 1/g$ around the zero. Then, for $w,z,\alpha,\beta \ll 1/g$, we approximate the integrand on these contours via the Taylor expansion

\begin{align}
    & \frac{(w+z) \zeta_q(1+w+z) A(\alpha,\beta;w,z)}{\zeta_q(1+\alpha+w) \zeta_q(1+\alpha+z) \zeta_q(1+\beta+w) \zeta_q(1+\beta+z)} \nonumber \\
    &\qquad =(\alpha+w) (\alpha+z) (\beta+w) (\beta+z) A(0,0;0,0) \log^3 q+O(g^{-5}).
\end{align}
Since $A(0,0;0,0)=1$ as can be seen by the definition, this gives us

\begin{align} \label{J(a,b) int}
    J_{\alpha,\beta}(y) & =\frac{\log^3 q}{(2\pi i)^2} \int_1^y \oint \oint \frac{u^{w+z}}{w^{m+1} z^{n+1}} \left( (\alpha+w) (\alpha+z) (\beta+w) (\beta+z)+O(g^{-5}) \right) \, dw \, dz \frac{du}{u} \nonumber \\
    & =\frac{\log^3 q}{(2\pi i)^2} \int_1^y \bigg( \oint \frac{u^w}{w^{m+1}} (\alpha+w) (\beta+w) \, dw \bigg) \bigg( \oint \frac{u^z}{z^{n+1}} (\alpha+z) (\beta+z) \, dz \bigg) \frac{du}{u} +O(g^{m+n-4}).
\end{align}
By generalising \eqref{residue 1} and \eqref{residue 2}, for $i \in \{0, 1\}$ and $j \in \{0, 1, 2\}$, we have that

\begin{equation} \label{general residue}
    \sum_{n \geq 1} \frac{p_{i,n} n!}{\log^n y} \frac{1}{2\pi i} \oint \frac{y^z}{z^{n+1-j}} \, dz=\frac{1}{\log^j y} P_i^{(j)}(1).
\end{equation}
Hence, by combining \eqref{I(a,b) def}, \eqref{J(a,b) int} and \eqref{general residue}, we have

\begin{align}
    I(\alpha,\beta) & =\zeta_q(1+2\alpha) \zeta_q(1+2\beta) \zeta_q(1+\alpha+\beta) \log^3 q \nonumber \\
    &\quad \times \Bigg[ \int_1^y \left(\alpha\beta P_1\left(\frac{\log u}{\log y}\right)+\frac{\alpha+\beta}{\log y} P_1'\left(\frac{\log u}{\log y}\right)+\frac{1}{\log^2 y} P_1''\left(\frac{\log u}{\log y}\right)\right) \nonumber \\
    &\qquad \times\left(\alpha\beta P_2\left(\frac{\log u}{\log y}\right)+\frac{\alpha+\beta}{\log y} P_2'\left(\frac{\log u}{\log y}\right)+\frac{1}{\log^2 y} P_2''\left(\frac{\log u}{\log y}\right)\right) \frac{du}{u}+O(g^{-4}) \Bigg].
\end{align}
By making the change of variables $u=y^r$, the integral in the above becomes

\begin{align}
    & \log y \int_0^1 \left(\alpha\beta P_1(r)+\frac{\alpha+\beta}{\log y} P_1'(r)+\frac{1}{\log^2 y} P_1''(r)\right) \left(\alpha\beta P_2(r)+\frac{\alpha+\beta}{\log y} P_2'(r)+\frac{1}{\log^2 y} P_2''(r)\right) dr.
\end{align}

Now, for $\alpha,\beta\asymp 1/g$ and $|\alpha+\beta|\gg 1/g$, we have the Laurent expansion

\begin{equation}
    \zeta_q(1+2\alpha) \zeta_q(1+2\beta) \zeta_q(1+\alpha+\beta)=\frac{1}{4\alpha\beta (\alpha+\beta) \log^3 q}+O(g^2)= O(g^3).
\end{equation}
Hence, we have that

\begin{align}
    I(\alpha,\beta) & =\frac{\log y}{4\alpha\beta (\alpha+\beta)} \int_0^1 \left(\alpha\beta P_1(r)+\frac{\alpha+\beta}{\log y} P_1'(r)+\frac{1}{\log^2 y} P_1''(r)\right) \nonumber \\ 
    &\qquad \times\left(\alpha\beta P_2(r)+\frac{\alpha+\beta}{\log y} P_2'(r)+\frac{1}{\log^2 y} P_2''(r)\right) dr+O(g^{-1}),
\end{align}
uniformly on any fixed annuli such that $\alpha,\beta\asymp 1/g$ and $|\alpha+\beta|\gg 1/g$. Returning to \eqref{M(a,b) decomp} with this formula and collecting similar terms gives us

\begin{align}
    \mathcal{M}(\alpha,\beta;P_1,P_2)& =\frac{1}{4} \left(\frac{\alpha\beta (1-q^{-2g(\alpha+\beta)})}{\alpha+\beta}+\frac{\alpha\beta (q^{-2g\alpha}-q^{-2g\beta})}{\alpha-\beta}\right) \log y \int_0^1 P_1(r) P_2(r) dr \nonumber \\
    &\quad +\frac{1}{4} (1+q^{-2g\alpha}) (1+q^{-2g\beta}) \int_0^1 \big(P_1(r) P_2'(r) +P_1'(r) P_2(r)\big) dr \nonumber \\
    &\quad +\frac{1}{4} \left(\frac{(1+q^{-2g\alpha}) (1-q^{-2g\beta})}{\beta}+\frac{(1-q^{-2g\alpha}) (1+q^{-2g\beta})}{\alpha}\right) \frac{1}{\log y} \int_0^1 P_1'(r) P_2'(r) dr \nonumber \\
    &\quad +\frac{1}{4} \left(\frac{1-q^{-2g(\alpha+\beta)}}{\alpha+\beta}+\frac{q^{-2g\beta}-q^{-2g\alpha}}{\alpha-\beta}\right) \frac{1}{\log y} \int_0^1 \big(P_1(r) P_2''(r)+P_1''(r) P_2(r)\big) dr \nonumber \\
    &\quad +\frac{(1-q^{-2g\alpha}) (1-q^{-2g\beta})}{4\alpha\beta} \frac{1}{\log^2 y} \int_0^1 \big(P_1'(r) P_2''(r)+P_1''(r) P_2'(r)\big) dr \nonumber \\
    &\quad +\frac{1}{4\alpha\beta} \left(\frac{1-q^{-2g(\alpha+\beta)}}{\alpha+\beta}+\frac{q^{-2g\alpha}-q^{-2g\beta}}{\alpha-\beta}\right) \frac{1}{\log^3 y} \int_0^1 P_1''(r) P_2''(r) dr+O(g^{-1}).
\end{align}
Finally, since both $\mathcal{M}(\alpha,\beta,P_1,P_2)$ and the main term on the right-hand side above are holomorphic for $\alpha,\beta\ll 1/g$, the error term is also holomorphic in this region. By the maximum modulus principle, the above formula holds uniformly for $\alpha,\beta\ll 1/g$ which completes the proof of the lemma.
\end{proof}

\subsection{Proof of Theorem \ref{2nd mollified moment thm}}

We define

\begin{equation}
    \mathcal{N}(\alpha,\beta;P_1,P_2)=\frac{1}{|\H|} \sum_{D \in \H} \Lambda(\tfrac{1}{2}+\alpha,\chi_D) \Lambda(\tfrac{1}{2}+\beta,\chi_D) M_y(\chi_D,P_1) M_y(\chi_D,P_2).
\end{equation}
By the definition of the completed $L$-function and Lemma \ref{2nd shifted mollified lemma}, we have

\begin{align}
    & \mathcal{N}(\alpha,\beta;P_1,P_2) \nonumber \\
    & =\frac{\log y}{4} \left(\frac{\alpha\beta (q^{g(\alpha+\beta)}-q^{-g(\alpha+\beta)})}{\alpha+\beta}+\frac{\alpha\beta (q^{g(\beta-\alpha)}-q^{g(\alpha-\beta)})}{\alpha-\beta}\right) \int_0^1 P_1(r) P_2(r) dr \nonumber \\
    &\quad +\frac{1}{4} (q^{g\alpha}+q^{-g\alpha}) (q^{g\beta}+q^{-g\beta}) \int_0^1 \big(P_1(r) P_2'(r) +P_1'(r) P_2(r)\big) dr \nonumber \\
    &\quad \frac{1}{4 \log y} \left(\frac{(q^{g\alpha}+q^{-g\alpha}) (q^{g\beta}-q^{-g\beta})}{\beta}+\frac{(q^{g\beta}-q^{-g\alpha}) (q^{g\beta}+q^{-g\beta})}{\alpha}\right) \int_0^1 P_1'(r) P_2'(r) dr \nonumber \\
    &\quad +\frac{1}{4 \log y} \left(\frac{q^{g(\alpha+\beta)}-q^{-g(\alpha+\beta)}}{\alpha+\beta}+\frac{q^{g(\alpha-\beta)}-q^{g(\beta-\alpha)}}{\alpha-\beta}\right) \int_0^1 \big(P_1(r) P_2''(r)+P_1''(r) P_2(r)\big) dr \nonumber \\
    &\quad +\frac{(q^{g\alpha}-q^{-g\alpha}) (q^{g\beta}-q^{-g\beta})}{4\alpha\beta} \frac{1}{\log^2 y} \int_0^1 \big(P_1'(r) P_2''(r)+P_1''(r) P_2'(r)\big) dr \nonumber \\
    &\quad +\frac{1}{4\alpha\beta} \left(\frac{q^{g(\alpha+\beta)}-q^{-g(\alpha+\beta)}}{\alpha+\beta}+\frac{q^{g(\beta-\alpha)}-q^{g(\alpha-\beta)}}{\alpha-\beta}\right) \frac{1}{\log^3 y} \int_0^1 P_1''(r) P_2''(r) dr+O(g^{-1}),
\end{align}
uniformly for $\alpha,\beta \ll 1/g$. Next, we scale the variables by writing $\alpha=\frac{a}{g \log q}$ and $\beta=\frac{b}{g \log q}$. Then, as $y=q^{2g\theta}$, the above can be rewritten as

\begin{align}
    \mathcal{N}(\alpha,\beta;P_1,P_2) & =\frac{1}{8\theta^3} \int_0^1 \frac{\sinh au}{a} \frac{\sinh bu}{b} du \int_0^1 P_1''(r) P_2''(r) dr \nonumber \\
    &\qquad +\frac{1}{4\theta^2} \frac{\sinh a}{a} \frac{\sinh b}{b} \int_0^1 \big(P_1'(r) P_2''(r)+P_1''(r) P_2'(r)\big)dr \nonumber \\
    &\qquad +\frac{1}{2\theta} \int_0^1 \cosh au \cosh bu \, du \int_0^1 \big(P_1(r) P_2''(r)+P_1''(r) P_2(r)\big) dr \nonumber \\
    &\qquad +\frac{1}{2\theta} \left(\frac{\sinh a \cosh b}{a}+\frac{\sinh b \cosh a}{b}\right) \int_0^1 P_1'(r) P_2'(r) dr \nonumber \\
    &\qquad +\cosh a \cosh b \int_0^1 \big(P_1(r) P_2'(r)+P_1'(r) P_2(r)\big) dr \nonumber \\
    &\qquad +2\theta ab \int_0^1 \sinh au \sinh bu \, du \int_0^1 P_1(r) P_2(r) dr+O(g^{-1}).
\end{align}
We now also define

\begin{equation}
    \mathcal{N}(Q_1,Q_2;P_1,P_2):=Q_1\left( \frac{d}{da} \right) Q_2\left( \frac{d}{db} \right) \left. \mathcal{N}\left( \frac{a}{g \log q},\frac{b}{g \log q},P_1,P_2 \right) \right|_{a=b=0}.
\end{equation}
For an even polynomial $Q$, we have

\begin{equation}
    Q \left( \frac{d}{da} \right) \left. \cosh{a} \right|_{a=0}=Q(1),
\end{equation}

\begin{equation}
    Q \left( \frac{d}{da} \right) \left. \frac{\sinh{au}}{a} \right|_{a=0}=Q \left( \frac{d}{da} \right) \left. \int_0^u \cosh{at} \, dt \right|_{a=0}=\int_0^u Q(t) \, dt,
\end{equation}
and

\begin{equation}
    Q \left( \frac{d}{da} \right) \left. a \sinh{au} \right|_{a=0}= Q \left( \frac{d}{da} \right) \left. \frac{d}{dt} \cosh{at} \right|_{a=0}=\frac{d}{dt} Q(t)=Q'(t).
\end{equation}
Using the above formulae, we then have

\begin{align} \label{N(P,Q)}
    \mathcal{N}(Q_1,Q_2;P_1,P_2) & =\frac{1}{8\theta^3} \int_0^1 \Tilde{Q}_1(u) \Tilde{Q}_2(u) du \int_0^1 P_1''(r) P_2''(r) dr \nonumber \\
    &\qquad +\frac{1}{4\theta^2} \Tilde{Q}_1(1) \Tilde{Q}_2(1) \int_0^1 \big(P_1'(r) P_2''(r)+P_1''(r) P_2'(r)\big)dr \nonumber \\
    &\qquad +\frac{1}{2\theta} \int_0^1 Q_1(u) Q_2(u) du \int_0^1 \big(P_1(r) P_2''(r)+P_1''(r) P_2(r)\big) dr \nonumber \\
    &\qquad +\frac{1}{2\theta} \left(Q_1(1) \Tilde{Q}_2(1)+\Tilde{Q}_1(1) Q_2(1)\right) \int_0^1 P_1'(r) P_2'(r) dr \nonumber \\
    &\qquad +Q_1(1) Q_2(1) \int_0^1 \big(P_1(r) P_2'(r)+P_1'(r) P_2(r)\big) dr \nonumber \\
    &\qquad +2\theta \int_0^1 Q_1'(u) Q_2'(u) du \int_0^1 P_1(r) P_2(r) dr+O(g^{-1}),
\end{align}
where

\begin{equation}
    \Tilde{Q}(u)=\int_0^u Q(t) \, dt.
\end{equation}

To complete the proof we write the above expression for $\mathcal{N}(P_1,P_2;Q_1,Q_2)$ in the more compact form given in the statement of Theorem \ref{2nd mollified moment thm}. This requires the following identities which all follow from integration by parts:

\begin{equation}
    \int_0^1 \big(P_1'(r) P_2(r)+P_1(r) P_2'(r) \big) \, dr=P_1(1) P_2(1),
\end{equation}

\begin{equation}
    \int_0^1 \big(P_1''(r) P_2'(r)+P_1'(r) P_2''(r) \big) \, dr=P_1'(1) P_2'(1),
\end{equation}

\begin{equation}
    \int_0^1 P_1''(r) P_2(r) \, dr=P_1'(1) P_2(1)-\int_0^1 P_1'(r) P_2'(r) \, dr,
\end{equation}
and

\begin{equation}
    \int_0^1 \Tilde{Q}_1(u) Q_2'(u) \, du=\Tilde{Q}(1) Q_2(1)-\int_0^1 Q_1(u) Q_2(u) \, du.
\end{equation}
where the last identity uses the fact that $\Tilde{Q}(0)=0$. Using these identities, it follows that

\begin{align}
    \mathcal{N}(Q_1,Q_2;P_1;P_2) & =2\theta \int_0^1 \int_0^1 P_1(r) Q_1'(u) P_2(r) Q_2'(u) \, du \, dr+P_1(1) P_2(1) Q_1(1) Q_2(1) \nonumber \\
    &\qquad +\frac{1}{2\theta} \big( P_1(1) Q_1(1) P_2'(1) \Tilde{Q}_2(1)+P_1'(1) \Tilde{Q}_1(1) P_2(1) Q_2(1) \big) \nonumber \\
    &\qquad -\frac{1}{2\theta} \Big( \int_0^1 \int_0^1 \big( P_1(r) Q_1'(u) P_2''(r) \Tilde{Q}_2(u)+P_1''(r) \Tilde{Q}_1(u) P_2(r) Q_2'(u) \big) \, du \, dr \Big) \nonumber \\
    &\qquad +\frac{1}{4\theta^2} \Tilde{Q}_1(1) \Tilde{Q}_2(1) P_1'(1) P_2'(1)+\frac{1}{8\theta^3} \int_0^1 \int_0^1 P_1''(r) \Tilde{Q}_1(u) P_2''(r) \Tilde{Q}_2(u) \, du \, dr \nonumber \\
    &\qquad +O(g^{-1}).
\end{align}
The expression on the right may then be factorised which yields the statement of Theorem \ref{2nd mollified moment thm}.

\section{Non-vanishing of $\Lambda^{(2k)}(\frac{1}{2},\chi_D)$}

In this section we will prove Theorem \ref{non-vanishing thm}. The strategy is to take $Q(x)=x^{2k}$ in Theorem \ref{general non-vanishing thm} and then find the optimal polynomial $P=P_{2k}$ to minimise the ratio $\mathcal{R}(P,x^{2k})$. For any polynomial $Q$, we denote

\begin{equation}
    I(Q)=\int_0^1 Q(t) \, dt.
\end{equation}
Recalling the expression for the ratio $\mathcal{R}(P,Q)$ given in Theorem \ref{non-vanishing thm}, we see that we need to minimise

\begin{equation}
    \mathcal{R}(P,Q)=\frac{\int_0^1 \big( (2\theta)^{-1} I(\Tilde{Q}^2) P''(r)^2-4\theta I(\Tilde{Q} Q') P(r) P''(r)+8\theta^3 I(Q'^2) P(r)^2 \big) \, dr}{\big( 2\theta P(1) Q(1)+P'(1) \Tilde{Q}(1) \big)^2},
\end{equation}
with $Q(x)$ ultimately taken to be $x^{2k}$. As mentioned in Remark \ref{non-vanishing remark}, the expression for the ratio $\mathcal{R}(P,Q)$ in \eqref{ratio def} is similar to the expression for the analogous ratio $\mathscr{R}_2(P,Q)$ considered in \cite{Kowalski2000}. Specifically, \cite[Eq. 32]{Kowalski2000} gives the expression

\begin{equation}
    \mathscr{R}_2(P,Q)=\frac{\int_0^1 \big( \Delta^{-1} I(Q^2) P''(r)^2-2\Delta I(Q Q'') P(r) P''(r)+\Delta^3 I(Q''^2) P(r)^2 \big) \, dr}{\big( \Delta P(1) Q'(1)+P'(1) Q(1) \big)^2}.
\end{equation}
Thus, by setting $\Delta=2\theta$, we have the relation $\mathcal{R}(P,Q)=\mathscr{R}_2(P,\Tilde{Q})$ and we may easily translate the results of the optimisation process carried out in \cite[Section 7]{Kowalski2000} to minimise $\mathscr{R}_2(P,Q)$ to our case considered here. For instance, in the case $k=0$ so $Q(x)=1$, we have $\Tilde{Q}(x)=x$ and it is shown in \cite{Kowalski2000} that the optimal $P$ which minimises $\mathscr{R}_2(P,x)$ is $P_0(x)=x^2-x^3/6$. We then compute that $\mathcal{R}(P_0,1)=1/7$ which by recalling \eqref{general non-vanishing bound}, gives us the lower bound $p_0=7/8$ on the proportion of $L(\tfrac{1}{2},\chi_D) \neq 0$. We note that this optimal polynomial and non-vanishing proportion is completely analogous to that found by Soundararajan in \cite{Soundararajan2000}.

For the case of general $k \geq 1$, the unique (up to scaling by a constant) choice of $P$ which minimises $\mathcal{R}(P,Q)=\mathscr{R}_2(P,\Tilde{Q})$ is of the form

\begin{equation} \label{optimal P}
    P(x)=\sinh(\alpha x) \big( \cos(\beta x)-Y\sin(\beta x) \big)-\frac{\alpha}{\beta} e^{-\alpha x} \sin(\beta x),
\end{equation}
where

\begin{align} \label{optimise alpha}
    \alpha=\frac{2\theta}{\sqrt{2 I(\Tilde{Q}^2)}} \sqrt{\sqrt{I(\Tilde{Q}^2) I(Q'^2)}+I(\Tilde{Q} Q')},
\end{align}
and

\begin{equation} \label{optimise beta}
    \beta=\frac{2\theta}{\sqrt{2 I(\Tilde{Q}^2)}} \sqrt{\sqrt{I(\Tilde{Q}^2) I(Q'^2)}-I(\Tilde{Q} Q')}.
\end{equation}
We note that $(\alpha \pm i\beta)^2$ are the roots of the polynomial

\begin{equation} \label{optimise 2}
    \frac{1}{2\theta} I(\Tilde{Q}^2) X^2-4\theta I(\Tilde{Q} Q') X+8 \theta^3 I(Q'^2).
\end{equation}
The value of parameter $Y$ is determined by the condition that $P(x)$ must satisfy

\begin{align} \label{optimise main}
    2\theta f(1) Q(1)+f'(1) \Tilde{Q}(1) & =f(1) \big( 8\theta^3 I(Q'^2) \Pi'(1)-2 \theta I(\Tilde{Q} Q') \Pi'''(1) \big) \nonumber \\
    &\quad +f'(1) \big( 2\theta I(\Tilde{Q} Q') \Pi''(1)-8\theta^3 I(Q'^2) \Pi(1) \big) \nonumber \\
    &\quad +\int_0^1 f''(r) \big( (2\theta)^{-1} I(\Tilde{Q}^2) \Pi''''(r)-4\theta I(\Tilde{Q} Q') \Pi''(r)+8 \theta^3 I(Q'^2) \Pi(r) \big) \, dr.
\end{align}
for all admissible functions $f$, where $\Pi(x)$ is a function with absolutely convergent Taylor series on $[0,1]$ such that $\Pi''(x)=P(x)$. One solves for $Y$ by substituting the optimal $P(x)$ in \eqref{optimal P} into \eqref{optimise main}.

With $Q(x)=x^{2k}$, we have

\begin{equation}
    \alpha=2\theta \sqrt{k(2k+1) \left( \sqrt{\frac{4k+3}{4k-1}}+\frac{4k+3}{4k+1} \right)},
\end{equation}
and

\begin{equation}
    \beta=2\theta \sqrt{k(2k+1) \left( \sqrt{\frac{4k+3}{4k-1}}-\frac{4k+3}{4k+1} \right)}.
\end{equation}
Also, the left-hand side of \eqref{optimise main} becomes

\begin{equation}
    2\theta f(1)+\frac{1}{2k+1} f'(1).
\end{equation}
So, we may substitute our optimal $P(x)$ above back into \eqref{optimise main} and set the ratio of the $f(1)$ and $f'(1)$ coefficients on the right-hand side to be $2\theta (2k+1)$ as well. This gives us a single linear equation in $Y$ and therefore a unique solution. While we do obtain an explicit formula for $Y$ for any $k \geq 1$ and $\theta>0$ using Mathematica, the formula is fairly messy and so we won't reproduce it here. Numerically for $\theta=1/2$, we compute the values of $Y$ for $k=1,2,3,4$ as $-0.8827$, $-0.7078$, $-0.6537$ and $-0.6268$. Evaluating $\mathcal{R}(P,x^{2k})$ using these values leads to the non-vanishing proportions given in Theorem \ref{non-vanishing thm}.

As mentioned above, the explicit formula for $Y$ is not concise. However, we do obtain a very simple formula for the asymptotic behaviour of $Y$ as $k \to \infty$. Using Mathematica, we have that as $k \to \infty$, $\alpha=2\theta (2k+1)+O(k^{-1})$ and $\beta=\theta+O(k^{-1})$, and we find that $Y$ is given asymptotically for large $k$ by

\begin{align}
    Y & =\frac{2 \alpha \beta \left(2 \alpha^2+\beta^2\right) \cos
   \beta-\left(4 \alpha^4+3 \alpha^2 \beta^2+\beta^4\right) \sin
   \beta}{2 \alpha \beta \left(2 \alpha^2+\beta^2\right) \sin
   \beta+\left(4 \alpha^4+3 \alpha^2 \beta^2+\beta^4\right) \cos \beta}+O(k^{-1}) \nonumber \\
    & =-\tan \beta+O(k^{-1}) \nonumber \\
    & =-\tan \theta+O(k^{-1}).
\end{align}

With the optimal choice for $P(x)$ now determined, we are ready to analyse the rate at which $\mathcal{R}(P,x^{2k})$ tends to zero.

\begin{lemma} \label{R decay}
For $P(x)=P_{2k}(x)$ the optimal choice of polynomial defined above and for any $0<\theta \leq 1/2$, we have that

    \begin{equation}
        \mathcal{R}(P,x^{2k})=\frac{1}{64 k^2}+O(k^{-3}),
    \end{equation}
    as $k \to \infty$.
\end{lemma}

\begin{proof}
The proof follows the same lines as that of \cite[Lemma 7.1]{Kowalski2000}. Recall that $\mathcal{R}(P,Q)$ is given by

\begin{equation}
    \mathcal{R}(P,Q)=\frac{\int_0^1 \big( (2\theta)^{-1} I(\Tilde{Q}^2) P''(x)^2-4\theta I(\Tilde{Q} Q') P(x) P''(x)+8\theta^3 I(Q'^2) P(x)^2 \big) \, dx}{\big( 2\theta P(1) Q(1)+P'(1) \Tilde{Q}(1) \big)^2}.
\end{equation}
With $Q(x)=x^{2k}$, the denominator of $\mathcal{R}(P,x^{2k})$ is the square of

\begin{equation}
    2\theta P(1)+\frac{1}{2k+1} P'(1).
\end{equation}
Using the expression for $P(x)$ in \eqref{optimal P}, we have that

\begin{align} \label{optimise denominator}
    2\theta P(1)+\frac{1}{2k+1} P'(1) & =\frac{e^{\alpha}}{2} \left( \left( 2\theta+\frac{\alpha}{(2k+1)} \right) (\cos\beta-Y \sin\beta) \right)+O(\alpha e^{-\alpha})+O(k^{-1} e^{\alpha}) \nonumber \\
    & =\frac{2e^{\alpha} \theta}{\cos\beta}+O(k^{-1} e^{\alpha}),
\end{align}
where we have used the estimates $\alpha=2\theta (2k+1)+O(k^{-1})$ and $Y=-\tan\beta+O(k^{-1})$ so that

\begin{equation}
    \cos\beta-Y \sin\beta=\frac{1}{\cos\beta}+O(k^{-1}).
\end{equation}

Now we consider the numerator of $\mathcal{R}(P,x^{2k})$. As $P(x)$ is given by a linear combination of the terms $e^{(\alpha+i\beta)x}$, the main term of the numerator will come from the $e^{+\alpha x}$ terms of $P(x)$ since the $e^{-\alpha x}$ terms and their derivatives, evaluated at $x=1$, will be exponentially small. Thus we may disregard the $e^{-\alpha x}$ terms of $P(x)$ when analysing the numerator of $\mathcal{R}(P,x^{2k})$ and therefore focus only on the contribution of

\begin{align}
    P_+(x) & :=\frac{e^{\alpha x}}{2} \big( \cos\beta-Y \sin\beta \big) \nonumber \\
    & =\frac{1}{4} \big( (1+iY) e^{(\alpha+i\beta)x}+(1-iY) e^{(\alpha-i\beta)x} \big)
\end{align}
to the numerator. Specifically, we are left to evaluate

\begin{equation} \label{optimise numerator}
    \int_0^1 \big( (2\theta)^{-1} I(\Tilde{Q}^2) P_+''(x)^2-4\theta I(\Tilde{Q} Q') P_+(x) P_+''(x)+8\theta^3 I(Q'^2) P_+(x)^2 \big) \, dx.
\end{equation}
As $(\alpha \pm i\beta)^2$ are the roots of the polynomial \eqref{optimise 2}, if we plug one of the individual terms of $P_+(x)$ into \eqref{optimise numerator} we get zero. Hence the terms which give a non-zero contribution to \eqref{optimise numerator} are those which involve a product of the two terms of $P_+(x)$. The contribution of these terms is

\begin{align} \label{optimise 3}
    \frac{1+Y^2}{16} & \Big( 2 (2\theta)^{-1} I(\Tilde{Q}^2) (\alpha^2+\beta^2)^2-4\theta I(\Tilde{Q} Q') \big( (\alpha+i\beta)^2+(\alpha-i\beta)^2 \big)+16\theta^3 I(Q'^2) \Big) \int_0^1 e^{2\alpha x} \, dx.
\end{align}
Also as $(\alpha \pm i\beta)^2$ are the roots of the polynomial \eqref{optimise 2}, we have that

\begin{equation}
    4\theta I(\Tilde{Q} Q') (\alpha \pm i\beta)^2=\frac{1}{2\theta} I(\Tilde{Q}^2) (\alpha \pm i\beta)^4+8\theta^3 I(Q'^2),
\end{equation}
and substituting this into \eqref{optimise 3} gives us

\begin{equation}
    \frac{1+Y^2}{16} \frac{I(\Tilde{Q}^2)}{2\theta} \left( 2(\alpha^2+\beta^2)^2-(\alpha+i\beta)^4-(\alpha-i\beta)^4 \right) \int_0^1 e^{2\alpha x} \, dx.
\end{equation}
Evaluating the integral and then simplifying this expression then yields

\begin{equation}
    \frac{e^{2\alpha}}{2\alpha} \frac{1+Y^2}{2\theta} I(\Tilde{Q}^2) (\alpha \beta)^2 \left( 1+O(e^{-2\alpha}) \right).
\end{equation}
Lastly, using $\alpha=2\theta (2k+1)+O(k^{-1})$, $\beta=\theta+O(k^{-1}$), $1+Y^2=(\cos\beta)^{-2}+O(k^{-1})$ and plugging in

\begin{equation}
    I(\Tilde{Q}^2)=\frac{1}{(4k+3) (2k+1)^2},
\end{equation}
gives us the final expression

\begin{equation}
    \frac{e^{2\alpha} \theta^2}{2 \cos^2 \beta} \frac{1}{(4k+3) (2k+1)} \left( 1+O(k^{-1}) \right)=\frac{e^{2\alpha} \theta^2}{16 k^2 \cos^2 \beta} \left( 1+O(k^{-1}) \right)
\end{equation}
for the numerator of $\mathcal{R}(P,x^{2k})$. Recalling that the denominator of $\mathcal{R}(P,x^{2k})$ is given by the square of \eqref{optimise denominator}, we therefore have that

\begin{equation}
    \mathcal{R}(P,x^{2k})=\frac{1}{64 k^2} \left( 1+O(k^{-1}) \right),
\end{equation}
as required.
\end{proof}

The proof of Theorem \ref{non-vanishing thm} is completed by recalling the statement of Theorem \ref{general non-vanishing thm} with $Q(x)=x^{2k}$ and applying the result of Lemma \ref{R decay}.

\section*{Acknowledgements}

The first author is grateful to the Leverhulme Trust (RPG-2017-320) for the support through the research project grant ``Moments of $L$-functions in Function Fields and Random Matrix Theory". The research of the second author was supported by an EPSRC Standard Research Studentship (EP/V520317/1) at the University of Exeter.

\section*{Open Access}

For the purpose of open access, the authors have applied a Creative Commons Attribution (CC BY) licence to any Author Accepted Manuscript version arising from this submission.

\bibliographystyle{plain}
\bibliography{Bibliography}

\end{document}